\documentclass[12pt]{article}
\usepackage{amsfonts,amsmath,amssymb,amscd}
\usepackage{shadow}
\usepackage{graphicx}
\usepackage{graphicx}
\usepackage{color}

\newtheorem{theo}{Theorem}[section]

\newtheorem{lemm}[theo]{Lemma}
\newtheorem{prop}[theo]{Proposition}

\def\qed{\hfill \rule{4pt}{7pt}}

\def\proof{\noindent {\it{Proof.} \hskip 2pt}}

\begin{document}
\begin{center}
{\large\bf The Generating Function for the Dirichlet Series $L_m(s)$}
\end{center}

\begin{center}
William Y.C. Chen$^1$, Neil J.Y. Fan$^2$, Jeffrey Y.T. Jia$^3$

Center for Combinatorics, LPMC-TJKLC\\
Nankai University, Tianjin 300071, P.R. China

$^1$chen@nankai.edu.cn, $^2$fjy@cfc.nankai.edu.cn,
$^3$jyt@cfc.nankai.edu.cn.

\end{center}

\vskip 6mm \noindent{\small {\bf Abstract.} The Dirichlet series
$L_m(s)$ are of fundamental importance in number theory.  Shanks defined the generalized Euler and class numbers in connection with
these Dirichlet series, denoted by
$\{s_{m,n}\}_{n\geq 0}$. We obtain a formula for the exponential generating function $s_m(x)$ of $s_{m,n}$, where $m$ is an arbitrary
 positive integer. In particular, for $m>1$, say, $m=bu^2$, where $b$ is square-free and $u>1$, we prove that $s_m(x)$ can be expressed as a linear combination of the four
 functions $w(b,t)\sec (btx)(\pm \cos ((b-p)tx)\pm \sin (ptx))$,
  where  $p$ is an integer satisfying $0\leq p\leq b$, $t|u^2$ and $w(b,t)=K_bt/u$ with $K_b$ being a constant depending on $b$.
Moreover, the Dirichlet series  $L_m(s)$
 can be easily computed from the generating function formula
  for $s_m(x)$.  Finally, we show that the main ingredient in the
   formula for $s_{m,n}$ has a combinatorial interpretation in terms of   the $m$-signed permutations defined by Ehrenborg and Readdy. In principle, this  answers a question posed by Shanks concerning a combinatorial interpretation for the numbers $s_{m,n}$ .



\section{Introduction}

The Dirichlet series
\begin{equation}
L_m(s)=\sum\limits_{l>0\atop {\rm odd}\ l}\left(\frac{-m}{l}\right)\frac{1}{l^{s}},\label{Diri}
\end{equation}
where $(-m/l)$ is the Jacobi symbol, originate in the distribution of primes into arithmetic progressions, the class number of binary quadratic forms, as well as the distribution of Legendre and Jacobi symbols. They play a crucial role in the computation of certain number-theoretic constants, see \cite{Bateman, Hardy, Pall,   Shanks3}. Several approaches have been developed for the computation of $L_m(s)$, see, for example, Shanks \cite{Shanks2, Shanks5, Shanks4}.

The generalized Euler and class numbers were introduced by Shanks for the
computation of the Dirichlet series $L_m(s)$ \cite{Shanks2, Shanks}.
These numbers are also  related to derivative polynomials and Euler polynomials, see Hoffman \cite{Hoff} and Shanks \cite{Shanks}.

In this paper, we obtain  the generating functions for the generalized Euler and class numbers. Let us recall the definition of the generalized Euler and class numbers $s_{m,n}$ $(m\geq 1, n\geq 0)$, introduced by Shanks,
\[
 s_{m,n}=\left\{
       \begin{array}{ll}
          c_{m,\frac{n}{2}} &\mbox{ if $n$ is even,}\\[5pt]
          d_{m,\frac{n+1}{2}} &\mbox{ if $n$ is odd.}
       \end{array}
       \right.
\]
where $c_{m,n}$ and $d_{m,n}$ are given by
\begin{equation}
c_{m,n}=(2n)!L_m(2n+1)(K_m\sqrt m)^{-1}\left(\frac{\pi}{2m}\right)^{-2n-1},\label{c}
\end{equation}
\begin{equation}
d_{m,n}=(2n-1)!L_{-m}(2n)(K_m\sqrt m)^{-1}\left(\frac{\pi}{2m}\right)^{-2n},\label{d}
\end{equation}
in which $K_m=\frac{1}{2}$ if $m=1$ and $1$ otherwise, and
the Dirichlet series $L_m(s)$ are defined by (\ref{Diri}).

Set
\begin{align*}
c_{m}(x)&=\sum\limits_{n\geq 0}c_{m,n}\frac{x^{2n}}{(2n)!},\\
d_{m}(x)&=\sum\limits_{n\geq 1}d_{m,n}\frac{x^{2n-1}}{(2n-1)!},\\
s_{m}(x)&=\sum\limits_{n\geq 0}s_{m,n}\frac{x^n}{n!}.
\end{align*}
Clearly,
\begin{equation}
s_{m}(x)=c_{m}(x)+d_{m}(x).\nonumber
\end{equation}
By the definitions (\ref{c}) and (\ref{d}), we have
\begin{align*}
\sum\limits_{n\geq 0}L_m(2n+1)x^{2n}&=\frac{K_m\pi\sqrt m}{2m}c_m\left(\frac{\pi}{2m}x\right),\\
\sum\limits_{n\geq 1}L_{-m}(2n)x^{2n-1}&=\frac{K_m\pi\sqrt m}{2m}d_m\left(\frac{\pi}{2m}x\right).
\end{align*}
Therefore, if we set
\[
 \hat{L}_{m}(s)=\left\{
       \begin{array}{ll}
          L_{-m}(s+1) &\mbox{ if $s$ is odd;}\\[5pt]
          L_{m}(s+1) &\mbox{ if $s$ is even,}
       \end{array}
       \right.
\]
then
\begin{equation}
\sum\limits_{s\geq 0}\hat{L}_{m}(s)x^{s}=\frac{K_m\pi\sqrt m}{2m}s_m\left(\frac{\pi}{2m}x\right).\label{sl}
\end{equation}
It follows from (\ref{sl}) that  $\hat{L}_{m}(s)$ is determined by $s_m\left(x\right)$. In other words,
 the generating function $s_m(x)$ leads to a quick way
  to compute $\hat{L}_m(s)$.

Consider $c_{m,n}$ and $d_{m,n}$ as entries of the infinite matrices $C$ and $D$, respectively. Then the first column of $C$ forms the sequence of
    class numbers in connection with primitive binary quadratic forms,
    and the first row of $C$ forms the sequence of secant numbers,
    corresponding to up-down permutations of even length.
   Meanwhile, the first row of $D$ forms the sequence of tangent numbers corresponding
     to up-down permutations of odd length. Recall that both secant numbers and
      tangent numbers are called Euler numbers. This is why the numbers $s_{m,n}$
         are called generalized Euler and class numbers.

Shanks \cite{Shanks} found  recurrence relations for $c_{m,n}$ and $d_{m,n}$ over the index $n$, from which it follows that $c_{m,n}$ and $d_{m,n}$ are integers.
 For example, we have
\[
\sum\limits_{i=0}^{n}(-4)^{i}{2n \choose 2i}c_{2,n-i}=(-1)^n
\]
and
\[
\sum\limits_{i=0}^{n-1}(-4)^{i}{2n-1 \choose 2i}d_{2,n-i}=(-1)^{n-1}.
\]
In fact, due to the well-known Euler product of the Dirichlet series $L_m(s)$ (see \cite{Class,Nath}), it can be easily shown that $c_{m,n}$ and $d_{m,n}$ are
positive.

 Shanks \cite{Shanks} raised the following question:
 Whether all of the generalized Euler and class numbers may  have some combinatorial interpretation?
The combinatorial interpretations of $s_{m,n}$ for $m=1,2,3,4$ have been found. Let $(s_{m,n})_{n\geq 0}$ denote the sequence \[ (c_{m,0},d_{m,1},c_{m,1},d_{m,2},c_{m,2},d_{m,3},\ldots).\]

For $m=1$, the sequence $(1,1,1,2,5,16,\ldots)$ is listed as $A000111$ in the
datbase of Sloane \cite{Sloane},
  which is called the sequence of Euler numbers, enumerating the number of alternating permutations on $[n]=\{1, 2, \ldots, n\}$.

For $m=2$, the sequence  $(1,1,3,11,57,361,\ldots)$ is numbered $A001586$ in \cite{Sloane}, which is also called the sequence of
        Springer numbers which arise in the work of
         Springer on the theory of Weyl group.

For $m=3$, the sequence $(1,2,8,46,352,3362,\ldots)$ is referred to $A007289$ in \cite{Sloane},
               and we call it the sequence of Ehrenborg and Readdy
               numbers because of their discovery of a combinatorial interpretation in terms of  alternating $3$-signed permutations on $[n]$, see \cite{Ehren}.

For $m=4$,  a combinatorial interpretation of the sequence $(1,4,16,128,1280,16384,\ldots)$ has been given implicitly by  Ehrenborg and Readdy
\cite{Ehren2} in terms of  non-augmented Andr\'e $4$-signed permutations on $[n]$.

For $m\geq 5$, we compute the generating functions $s_{m}(x)$. For $m=1,2,3,4$, it is known that
\begin{align*}
s_1(x)&=\sec x+\tan x,\\
s_2(x)&=\frac{\cos x+\sin x}{\cos 2x},\\
s_3(x)&=\frac{\sin 2x+\cos x}{\cos 3x},\\
s_4(x)&=\sec 4x+\tan 4x.
\end{align*}
From our general formula, we get the following expressions for $m=5, 6, 7$,
\begin{align*}
s_5(x)&=\frac{\cos 4x+\sin x}{\cos 5x}+\frac{\cos 2x+\sin 3x}{\cos
5x},\\
s_6(x)&=\frac{\cos 5x+\sin x}{\cos 6x}+\frac{\cos x+\sin 5x}{\cos
6x},\\
s_7(x)&=\frac{\cos 3x+\sin 4x}{\cos 7x}+\frac{\cos x+\sin 6x}{\cos
7x}-\frac{\cos 5x+\sin 2x}{\cos 7x}.\nonumber
\end{align*}

Our paper is organized as follows. In Section \ref{sec2}, we compute the generating
 function $s_m(x)$ when $m$ is square-free, while in Section \ref{sec3}
 we consider the case when $m$ is not square-free. Section \ref{sec4} is devoted to the
   combinatorial interpretation of the numbers $s_{m,n}$ in terms of $m$-signed permutations as introduced by Ehrenborg and Readdy.

\section{Computation for $s_m(x)$ when $m>1$ is square-free}\label{sec2}

In this section, we compute the
generating function  $s_m(x)$ when $m>1$ is square-free. For $0\leq p\leq m$,
we follow the following notation used in \cite{Ehren},
\begin{equation}
\Lambda_{m,p}(x):=\frac{\cos ((m-p)x)+\sin (px)}{\cos (mx)}.\label{exponent}
\end{equation}

When $m$ is square-free, we shall not encounter the case that $m$ is a multiple
of $4$.  We shall have three formulas for  $s_m(x)$ depending on the
residue of $m$ modulo $4$.

\begin{theo}\label{3}
Assume that $m$ is square-free and  $m=4t+3$. Then
\begin{equation}
s_m(x)=\sum\limits_{k=1}^{t}\left(\frac{k}{m}\right)\Lambda_{m,4k}(x)+\sum_{k=t+1}^{2t+1}\left(\frac{k}{m}\right)\Lambda_{m,2m-4k}(x).\label{free3}
\end{equation}
\end{theo}

\begin{theo}\label{1}
Assume that $m$ is square-free and $m=4t+1$. Then
\begin{equation}
s_m(x)=\sum\limits_{k=1}^{t}\left(\frac{k}{m}\right)\Lambda_{m,m-4k}(x)-\sum_{k=t+1}^{2t}\left(\frac{k}{m}\right)\Lambda_{m,4k-m}(x).\label{free1}
\end{equation}
\end{theo}

\begin{theo}\label{2}
Assume that $m$ is square-free and $m=4t+2$. Then
\begin{equation}
s_m(x)=\sum\limits_{k=1\atop {\rm odd}\ k}^{4t+1}\left(\frac{-m}{k}\right)\Lambda_{m,k}(x).\label{free2}
\end{equation}
\end{theo}

To prove the above theorems, let us first recall the following formula of $L_m(2n+1)$ obtained by Shanks \cite{Shanks2, Shanks}.

\begin{lemm}\label{2n+1}
Suppose that $m>1$ is square-free. Then $L_m(2n+1)$ can be expressed as a linear combination of the Fourier series $S_{2n+1}(x)$ as follows
\begin{align*}
L_m(2n+1)=\frac{2}{\sqrt m}\sum\limits_k\epsilon_kS_{2n+1}(y_k),
\end{align*}
where the Jacobi symbols $\epsilon_k$ and rational numbers $y_k$ are uniquely determined by $m$, and $S_{2n+1}(x)$ is defined by
\[
S_{2n+1}(x)=\sum\limits_{k=0}^{\infty}\frac{\sin 2\pi(2k+1)x}{(2k+1)^{2n+1}}.
\]
Furthermore, we have
\begin{equation}
c_m(x)=\frac{1}{\cos (mx)}\sum\limits_k\epsilon_k\cos (mx(1-4y_k)).\nonumber
\end{equation}
\end{lemm}

In fact, Shanks has given an explicit procedure to determine the constants
$\epsilon_k$ and $y_k$. To compute $\epsilon_k$ and $y_k$, we use the definition (\ref{Diri}) of the series $L_m(s)$ and express the Jacobi symbol $\left(\frac{-m}{l}\right)$ as a linear combination of sines according to the following expansion, see \cite{Landau}.

\begin{prop}\label{expand}
Assume that $l$ is odd and $m$ satisfies the following two conditions: $m\equiv 1\ (\bmod\  4)$ or $m\equiv 8$ or $12\ (\bmod\  16)$ and $p^2\nmid m$ for any odd prime $p$. Then  we have
\begin{equation}
\left(\frac{m}{l}\right)=\frac{1}{\sqrt {m}}\sum\limits_{r=1}^{|m|}\left(\frac{m}{r}\right)e^{2\pi ilr/|m|}.\label{Landau}
\end{equation}
\end{prop}

In particular, when $m\equiv 3\ (\bmod\  4)$, then we have $-m\equiv 1\ (\bmod\  4)$ and  we can use the above expansion for $\left(\frac{-m}{l}\right)$. Similarly, when $m\equiv 1\ (\bmod\  4)$, we have  $-4m\equiv 12\ (\bmod\  16)$ so that we can compute $\left(\frac{-4m}{l}\right)$ by using the above formula. Finally, when $m\equiv 2\ (\bmod\  4)$, we have $-4m\equiv 8\ (\bmod\  16)$ so that we can compute $\left(\frac{-4m}{l}\right)$. Note that  when $l$ is odd, we have
\[
\left(\frac{-4m}{l}\right)=\left(\frac{-m}{l}\right).
\]
Thus, the Jacobi symbol $\left(\frac{-m}{l}\right)$  can be determined by the above procedure for $m> 1$.

On the other hand,  Shanks \cite{Shanks2, Shanks} provided the following formula for  $L_{-m}(2n)$.

\begin{lemm}\label{2n}
Suppose that $m>1$ is square-free. Then $L_{-m}(2n)$ can be expressed as a linear combination of the Fourier series $C_{2n}(x)$ as follows
\begin{align*}
L_{-m}(2n)=\frac{2}{\sqrt m}\sum\limits_k\epsilon_k'C_{2n}(y_k'),
\end{align*}
where the Jacobi symbols $\epsilon_k'$ and rational numbers  $y_k'$ are uniquely determined by $m$, and $C_{2n}(x)$ is defined by
\[
C_{2n}(x)=\sum\limits_{k=0}^{\infty}\frac{\cos 2\pi(2k+1)x}{(2k+1)^{2n}}.
\]
Furthermore, we have
\begin{equation}
d_m(x)=\frac{1}{\cos (mx)}\sum\limits_k\epsilon_k'\sin (mx(1-4y_k')). \nonumber
\end{equation}
\end{lemm}

Similarly, Shanks has shown how to compute the constants $\epsilon_k'$ and $y_k'$. In order to compute $\epsilon_k'$ and $y_k'$, we use the definition of the series $L_{-m}(s)$
\begin{equation}
L_{-m}(s)=\sum\limits_{l>0\atop {\rm odd}\ l}\left(\frac{m}{l}\right)\frac{1}{l^{s}},\nonumber
\end{equation}
and express the Jacobi symbol $\left(\frac{m}{l}\right)$ as a linear combination of cosines resorting  to  proposition \ref{expand}.

 For the case $m\equiv 3\ (\bmod\  4)$, we have $4m\equiv 12\ (\bmod\  16)$ so
 that  we can use the above expansion for $\left(\frac{4m}{l}\right)$.
 When $m\equiv 1\ (\bmod\  4)$, we can also compute $\left(\frac{m}{l}\right)$ by using the above expansion. Finally, when $m\equiv 2\ (\bmod\  4)$, we have $4m\equiv 8\ (\bmod\  16)$ and $\left(\frac{4m}{l}\right)$ can be determined
 in the same manner. Note that  when $l$ is odd, we have
\[
\left(\frac{4m}{l}\right)=\left(\frac{m}{l}\right).
\]
Thus, the Jacobi symbol $\left(\frac{m}{l}\right)$  can be determined for $m> 1$.

Keep in mind that $m$ is assumed to be square-free. Set
\[ \hat{c}_m(x)=\cos (mx)c_m(x), \quad \hat{d}_m(x)=\cos (mx)d_m(x), \quad \hat{s}_m(x)=\cos (mx)s_m(x).\]

\noindent {\it Proof of Theorem \ref{3}}. If $m\equiv 3\ (\bmod\  4)$, by using the expansion (\ref{Landau}) for $\left(\frac{-m}{l}\right)$ and $\left(\frac{4m}{l}\right)$, we have
\[\epsilon_k=\left(\frac{k}{m}\right), \quad y_k=\frac{k}{m}, \quad \epsilon_k'=\left(\frac{m}{k}\right),\quad y_k'=\frac{k}{4m},\]
which imply that
\begin{align*}
L_m(2n+1)&=\frac{2}{\sqrt m}\sum\limits_{k=1}^{(m-1)/2}\left(\frac{k}{m}\right)S_{2n+1}\left(\frac{k}{m}\right),\\
L_{-m}(2n)&=\frac{2}{\sqrt m}\sum\limits_{{\rm odd}\  k<m}\left(\frac{m}{k}\right)C_{2n}\left(\frac{k}{4m}\right).
\end{align*}
Therefore, by Lemmas \ref{2n+1} and \ref{2n}, we have
\begin{align*}
\hat{s}_m(x)=\hat{c}_m(x)+\hat{d}_m(x)
=\sum\limits_{k=1}^{(m-1)/2}
\left(\frac{k}{m}\right)\cos(m-4k)x+\sum\limits_{{\rm odd}\
k<m}\left(\frac{m}{k}\right)\sin(m-k)x.
\end{align*}
Suppose that $m=4t+3$. It follows that
\begin{align*}
\hat{c}_m(x)&=\sum\limits_{k=1}^{t}\left(\frac{k}{m}\right)\cos(m-4k)x+\sum\limits_{k=t+1}^{2t+1}\left(\frac{k}{m}\right)\cos(m-4k)x\\
&=\sum\limits_{k=1}^{t}\left(\frac{k}{m}\right)\cos(m-4k)x+\sum\limits_{k=t+1}^{2t+1}\left(\frac{k}{m}\right)\cos(4k-m)x,\\
\hat{d}_m(x)&=\sum\limits_{k=1 \atop {\rm odd} \ k
}^{4t+1}\left(\frac{m}{k}\right)\sin(m-k)x.
\end{align*}
Thus we obtain
\begin{align*}
\hat{s}_m(x)=&\sum\limits_{k=1}^{t}\left(\left(\frac{k}{m}\right)\cos(m-4k)x+\left(\frac{m}{m-4k}\right)\sin(4k)x\right)\\
&+\sum\limits_{k=t+1}^{2t+1}\left(\left(\frac{k}{m}\right)\cos(4k-m)x+\left(\frac{m}{4k-m}\right)\sin
(2m-4k)x\right).
\end{align*}
It remains  to verify that
\begin{equation}
\left(\frac{k}{m}\right)=\left(\frac{m}{m-4k}\right)\label{1t3}
\end{equation}
for $1\leq k\leq t$ and
\begin{equation}
\left(\frac{k}{m}\right)=\left(\frac{m}{4k-m}\right)\label{t+1}
\end{equation}
for $t+1\leq k\leq 2t+1$.

With respect to (\ref{1t3}),
since both $m$ and $m-4k$ are odd positive numbers, by the law of quadratic reciprocity, we find
\begin{align*}
\left(\frac{m}{m-4k}\right)&=-\left(\frac{m-4k}{m}\right)=-\left(\frac{-4k}{m}\right)=-\left(\frac{-4}{m}\right)\left(\frac{k}{m}\right)\\
&=-\left(\frac{-1}{m}\right)\left(\frac{4}{m}\right)\left(\frac{k}{m}\right)=\left(\frac{2}{m}\right)^2\left(\frac{k}{m}\right)=\left(\frac{k}{m}\right).
\end{align*}
Similarly,  (\ref{t+1}) can be checked via the following steps,
\begin{align*}
\left(\frac{m}{4k-m}\right)&=\left(\frac{4k-m}{m}\right)=\left(\frac{-(m-4k)}{m}\right)=\left(\frac{-1}{m}\right)\left(\frac{m-4k}{m}\right)\\
&=-\left(\frac{m-4k}{m}\right)=-\left(\frac{-4k}{m}\right)=\left(\frac{k}{m}\right).
\end{align*}
Hence we deduce that
\begin{align*}
s_m(x)=&\sum\limits_{k=1}^{t}\left(\frac{k}{m}\right)\frac{\cos \left(m-4k\right)x+\sin \left(4k\right)x}{\cos mx}\\
&+\sum\limits_{k=t+1}^{2t+1}\left(\frac{k}{m}\right)\frac{\cos(4k-m)x+\sin
(2m-4k)x}{\cos mx}.
\end{align*}
This completes the proof.\qed

\noindent {\it Proof of Theorem \ref{1}}. If $m\equiv 1\ (\bmod\  4)$, applying the  expansion (\ref{Landau}) to $\left(\frac{-4m}{l}\right)$ and $\left(\frac{m}{l}\right)$, we get
\[
\epsilon_k=\left(\frac{-m}{k}\right), \quad y_k=\frac{k}{4m}, \quad \epsilon_k'=\left(\frac{k}{m}\right),\quad y_k'=\frac{k}{m}.
\]
It follows that
\begin{align*}
L_m(2n+1)&=\frac{2}{\sqrt m}\sum\limits_{{\rm odd}\  k<m}\left(\frac{-m}{k}\right)S_{2n+1}\left(\frac{k}{4m}\right),\\
L_{-m}(2n)&=\frac{2}{\sqrt m}\sum\limits_{k=1}^{(m-1)/2}\left(\frac{k}{m}\right)C_{2n}\left(\frac{k}{m}\right).
\end{align*}
In view of Lemmas \ref{2n+1} and \ref{2n}, we find
\begin{equation}
\hat{s}_m(x)=\hat{c}_m(x)+\hat{d}_m(x)=\sum\limits_{{\rm odd}\
k<m}\left(\frac{-m}{k}\right)\cos(m-k)x+\sum\limits_{k=1}^{(m-1)/2}\left(\frac{k}{m}\right)\sin(m-4k)x.\nonumber
\end{equation}
Writing $m=4t+1$, we obtain
\begin{align*}
\hat{d}_m(x)&=\sum\limits_{k=1}^{t}\left(\frac{k}{m}\right)\sin(m-4k)x+\sum\limits_{k=t+1}^{2t}\left(\frac{k}{m}\right)\sin(m-4k)x\\
&=\sum\limits_{k=1}^{t}\left(\frac{k}{m}\right)\sin(m-4k)x-\sum\limits_{k=t+1}^{2t}\left(\frac{k}{m}\right)\sin(4k-m)x,\\
\hat{c}_m(x)&=\sum\limits_{k=1 \atop {\rm odd}\ k
}^{4t-1}\left(\frac{-m}{k}\right)\cos(m-k)x.
\end{align*}
It follows that
\begin{align*}
\hat{s}_m(x)=&\sum\limits_{k=1}^{t}\left(\left(\frac{k}{m}\right)\sin(m-4k)x+\left(-\frac{m}{m-4k}\right)\cos (4k)x\right)\\
&+\sum\limits_{k=t+1}^{2t}\left(-\left(\frac{k}{m}\right)\sin(4k-m)x+\left(-\frac{m}{4k-m}\right)\cos
(2m-4k)x\right).
\end{align*}
Finally, we need to show that
\begin{equation}
\left(\frac{k}{m}\right)=\left(-\frac{m}{m-4k}\right)\label{1t1}
\end{equation}
for $1\leq k\leq t$ and
\begin{equation}
-\left(\frac{k}{m}\right)=\left(-\frac{m}{4k-m}\right)\label{t}
\end{equation}
for $t+1\leq k\leq 2t$.

To confirm (\ref{1t1}),
since both $m$ and $m-4k$ are odd positive numbers, we may
 employ the law of quadratic reciprocity to deduce
\begin{align*}
\left(-\frac{m}{m-4k}\right)&=\left(\frac{-1}{m-4k}\right)\left(\frac{m}{m-4k}\right)=\left(\frac{m}{m-4k}\right)=\left(\frac{m-4k}{m}\right)=\left(\frac{-4k}{m}\right)\\
&=\left(\frac{-4}{m}\right)\left(\frac{k}{m}\right)=\left(\frac{-1}{m}\right)\left(\frac{2}{m}\right)^2\left(\frac{k}{m}\right)=\left(\frac{k}{m}\right).
\end{align*}
Moreover,  (\ref{t}) can be checked as follows
\begin{align*}
-\left(-\frac{m}{4k-m}\right)&=-\left(\frac{-1}{4k-m}\right)\left(\frac{m}{4k-m}\right)=\left(\frac{m}{4k-m}\right)=\left(\frac{4k-m}{m}\right)\\
&=\left(\frac{-(m-4k)}{m}\right)=\left(\frac{-1}{m}\right)\left(\frac{m-4k}{m}\right)
=\left(\frac{-4k}{m}\right)=\left(\frac{k}{m}\right).
\end{align*}
So we arrive at
\begin{align*}
s_m(x)=&\sum\limits_{k=1}^{t}\left(\frac{k}{m}\right)\frac{\cos(4k)x+\sin (m-4k)x}{\cos mx}\\
&-\sum\limits_{k=t+1}^{2t}\left(\frac{k}{m}\right)\frac{\cos(2m-4k)x+\sin
(4k-m)x}{\cos mx}.
\end{align*}
This completes the proof.\qed

\noindent {\it Proof of Theorem \ref{2}}. If $m\equiv 2\ (\bmod\  4)$, by using the expansion (\ref{Landau}) for $\left(\frac{-4m}{l}\right)$ and $\left(\frac{4m}{l}\right)$, we obtain
\[\epsilon_k=\left(\frac{-m}{k}\right), \quad y_k=\frac{k}{4m}, \quad \epsilon_k'=\left(\frac{m}{k}\right),\quad y_k'=\frac{k}{4m}.\]
Consequently,
\begin{align*}
L_m(2n+1)&=\frac{2}{\sqrt m}\sum\limits_{{\rm odd}\  k<m}\left(\frac{-m}{k}\right)S_{2n+1}\left(\frac{k}{4m}\right),\\
L_{-m}(2n)&=\frac{2}{\sqrt m}\sum\limits_{{\rm odd}\  k<m}\left(\frac{m}{k}\right)C_{2n}\left(\frac{k}{4m}\right).
\end{align*}
By Lemmas \ref{2n+1} and \ref{2n}, we see that
\begin{equation}
\hat{s}_m(x)=\hat{c}_m(x)+\hat{d}_m(x)=\sum\limits_{{\rm odd}\
k<m}\left(\frac{-m}{k}\right)\cos(m-k)x+\sum\limits_{{\rm odd}\
k<m}\left(\frac{m}{k}\right)\sin(m-k)x.\nonumber
\end{equation}
Writing $m=4t+2$, since
\[
\left(\frac{4t+2}{2t+1}\right)=0,
\]
we obtain
\begin{align*}
\hat{s}_m(x)=\sum\limits_{k=1 \atop {\rm odd}\ k
}^{4t+1}\left(\left(\frac{-m}{k}\right)\cos(m-k)x+\left(\frac{m}{m-k}\right)\sin
(k)x\right).
\end{align*}
Finally, it follows from $-k\equiv m-k\ (\bmod\  m)$ that
\begin{equation}
\left(\frac{-m}{k}\right)=\left(\frac{m}{m-k}\right)\label{14t+1}
\end{equation}
for $1\leq k\leq 4t+1$.
Hence we have reached the conclusion
\begin{align*}
s_m(x)=\sum\limits_{k=1 \atop {\rm odd} \ k
}^{4t+1}\left(\frac{-m}{k}\right)\frac{\cos(m-k)x+\sin(k)x}{\cos
mx}.
\end{align*}
This completes the proof.\qed

For $m=5,6,7$, the generating function $s_m(x)$ has been given in the
introduction.

\section{Computation for $s_m(x)$ when $m$ is not square-free}\label{sec3}

In this section, we obtain an expression for $s_m(x)$ for the case when
$m$ is not square-free. Assume that $m$ can be divided by a square $u^2>1$, say, $m=bu^2$,
where $b$ is square-free.
 Below  is the formula for this case.

\begin{theo}\label{lincombi}
Suppose that $m=bu^2$ as given above. Then we can express
$s_m(x)$ as a  linear combination of the four terms
\begin{equation}
w(b,t)\sec (btx)(\pm \cos ((b-p)tx)\pm \sin (ptx)),\nonumber
\end{equation}
where $p$ is an integer satisfying $0\leq p\leq b$  and $t|u^2$, and
the  coefficient $w(b,t)=K_bt/u$ is uniquely determined for any given $t$.
\end{theo}

The idea of the proof is to establish  two recursive relations (\ref{cbl}) and (\ref{dbl}) between $s_{m,n}$ and $s_{b,n}$.
Then we express $s_m(x)$ as a linear combination of the two terms $c_b(tx)$ and $d_b(tx)$ by considering the two cases according to
 whether there exist odd prime factors
 $u_i$ of $u$  with residues
$3$ modulo $4$. Since $b$ is square-free, $c_b(tx)$ and $d_b(tx)$ can be evaluated by using the formulas in the previous
section.

\proof Let us start with the following relation given by Shanks \cite{Shanks2}
\begin{equation}
L_m(s)=L_b(s)\prod\limits_{u_i|u}\left(1-\left(\frac{-b}{u_i}\right)\frac{1}{u_i^s}\right),\label{Facto}
\end{equation}
where the product ranges over  odd primes $u_i$ (if any) that divide $u$.
To be precise, in case there are no odd prime factors, the empty product is
assumed to be one by convention.
 From the definitions (\ref{c}) and (\ref{d})
  it follows that
\begin{align}
c_{m,n}&=K_bu(u^2)^{2n}\prod\limits_i\left(u_i^{2n+1}-\left(\frac{-b}{u_i}\right)\right)\left(\prod\limits_i\frac{1}{u_i}\right)^{2n+1}c_{b,n},
\label{cbl}\\
d_{m,n}&=K_bu(u^2)^{2n-1}\prod\limits_i\left(u_i^{2n}-\left(\frac{b}{u_i}\right)\right)\left(\prod\limits_i\frac{1}{u_i}\right)^{2n}d_{b,n}.\label{dbl}
\end{align}

For the purpose of computing $s_m(x)$ for the case when $m$ is
not square-free, we need to consider the two cases according to whether there exist $u_i\equiv 3\ (\bmod\  4)$
among the $k$ odd factors $u_1,u_2,\ldots,u_k$ of $u$.

Case 1: $u_i\equiv 1\ (\bmod\  4)$ for $1\leq i\leq k$.
In this case, we see that
\[
\left(\frac{-b}{u_i}\right)=\left(\frac{b}{u_i}\right).
\]
Suppose that among the factors $u_{1},u_{2},\ldots,u_{k}$
there are $k_1$ primes $u_1,u_2,\ldots,u_{k_1}$ satisfying
$\left(\frac{b}{u_i}\right)=1$ for $1\leq i\leq k_1$, $k_2$ primes
$u_{k_1+1},u_{k_1+2},\ldots,u_{k_1+k_2}$ satisfying
$\left(\frac{b}{u_{k_1+j}}\right)=-1$ for $1\leq j\leq k_2$,  and $k_3$
primes $u_{k_1+k_2+1},u_{k_1+k_2+2},\ldots,u_{k_1+k_2+k_3}$
satisfying $\left(\frac{b}{u_{k_1+k_2+l}}\right)=0$ for $1\leq l\leq
k_3$, where $k_1+k_2+k_3=k$. By (\ref{cbl}), we get
\begin{equation}
c_{m,n}=K_bu(u^{2})^{2n}\frac{\prod_{i=1}^{k_1}(u_{i}^{2n+1}-1)\times\prod_{j=1}^{k_2}(u_{k_1+j}^{2n+1}+1)}{\prod_{i=1}
^{k_1+k_2}u_{i}^{2n+1}}c_{b,n}.\label{cmn}
\end{equation}
Let
\begin{equation}
f_c=\frac{\prod_{i=1}^{k_1}(u_i^{2n+1}-1)\times
\prod_{j=1}^{k_2}(u_{k_1+j}^{2n+1}+1)}{\prod_{i=1}^{k_1+k_2}u_i^{2n+1}}.\nonumber
\end{equation}
In this notation, (\ref{cmn}) can be rewritten as
\begin{align*}
c_{m,n}=K_bu(u^{2})^{2n}f_cc_{b,n},
\end{align*}
which implies that
\begin{equation}
c_m(x)=K_bu\sum\limits_{n\geq 0}f_cc_{b,n}\frac{(u^2x)^{2n}}{(2n)!}.\label{cinsert}
\end{equation}
Since
\begin{equation}
\prod_{i=1}^{k_1}(u_i^{2n+1}-1)=\prod_{i=1}^{k_1}u_i^{2n+1}-\sum_{i=1}^{k_1}\left(\prod
u_1\cdots u_{i-1}u_{i+1}\cdots
u_{k_1}\right)^{2n+1}+\cdots+(-1)^{k_1}\nonumber
\end{equation}
and
\begin{equation}
\prod_{j=1}^{k_2}(u_{k_1+j}^{2n+1}+1)=\prod_{j=1}^{k_2}u_{k_1+j}^{2n+1}+\sum_{j=1}^{k_2}\left(\prod
u_{k_1+1}\ldots u_{k_1+j-1}u_{k_1+j+1}\cdots
u_{k_1+k_2}\right)^{2n+1}+\cdots+1,\nonumber
\end{equation}
we can expand $f_c$ as follows
\begin{align}
f_c=&1+\sum\limits_{j=1}^{k_2}\frac{1}{u_{k_1+j}^{2n+1}}+
\cdots+\frac{1}{\prod_{j=1}^{k_2}u_{k_1+j}^{2n+1}}
-\sum\limits_{i=1}^{k_1}\frac{1}{u_i^{2n+1}}
-\sum\limits_{i=1}^{k_1}
\sum\limits_{j=1}^{k_2}\frac{1}{u_i^{2n+1}u_{k_1+j}^{2n+1}}+\cdots\nonumber\\
&\;-\sum\limits_{i=1}^{k_1}\frac{1}{u_i^{2n+1}\prod_{j=1}^{k_2}u_{k_1+j}^{2n+1}}+\cdots+\frac{(-1)^{k_1}}{\prod_{i=1}^{k_1+k_2}u_i^{2n+1}}\label{fff}.
\end{align}

Plugging (\ref{fff}) into (\ref{cinsert}), we find that $c_m(x)$ is a linear
combination of the terms $c_b(tx)$, where $t|u^2$  and the  coefficient of $c_b(tx)$ in the linear combination is equal to  $K_bt/u$.

Similarly, we have
\begin{align}
d_{m,n}=K_bu(u^{2})^{2n-1}\frac{\prod_{i=1}^{k_1}(u_{i}^{2n}-1)\times\prod_{j=1}^{k_2}(u_{k_1+j}^{2n}+1)}
{\prod_{i=1}^{k_1+k_2}u_{i}^{2n}}d_{b,n}.\label{dmn}
\end{align}
Let
\begin{equation}
f_d=\frac{\prod_{i=1}^{k_1}(u_{i}^{2n}-1)\times\prod_{j=1}^{k_2}(u_{k_1+j}^{2n}+1)}
{\prod_{i=1}^{k_1+k_2}u_{i}^{2n}},\nonumber
\end{equation}
then  (\ref{dmn}) can be rewritten as
\begin{equation}
d_{m,n}=K_bu(u^{2})^{2n-1}f_dd_{b,n},\nonumber
\end{equation}
which leads to the following relation
\begin{equation}
d_m(x)=K_bu\sum\limits_{n\geq 1}f_dd_{b,n}\frac{(u^2x)^{2n-1}}{(2n-1)!}.\label{dinsert}
\end{equation}
Again we may expand $f_d$ as follows
\begin{align}
f_d=&1+\sum\limits_{j=1}^{k_2}\frac{1}{u_{k_1+j}^{2n}}+
\cdots+\frac{1}{\prod_{j=1}^{k_2}u_{k_1+j}^{2n}}
-\sum\limits_{i=1}^{k_1}\frac{1}{u_i^{2n}}
-\sum\limits_{i=1}^{k_1}
\sum\limits_{j=1}^{k_2}\frac{1}{u_i^{2n}u_{k_1+j}^{2n}}+\cdots\nonumber\\
&\;-\sum\limits_{i=1}^{k_1}\frac{1}{u_i^{2n}\prod_{j=1}^{k_2}u_{k_1+j}^{2n}}+\cdots+\frac{(-1)^{k_1}}
{\prod_{i=1}^{k_1+k_2}u_i^{2n}}\label{f'}.
\end{align}

Substituting (\ref{f'}) into (\ref{dinsert}), we find again  that $d_m(x)$ is a linear combination of the terms $d_b(tx)$, where $t|u^2$ the
 coefficient of $d_b(tx)$ in the linear combination is equal to  $K_bt/u$.
Furthermore, we see that $c_m(x)$ and $d_m(x)$ have the same
coefficients for the linear combinations. In other words, the relation
for $c_m(x)$ and $c_b(tx)$ is still valid after changing $c_m(x)$ and $c_b(tx)$
to $d_m(x)$ and $d_b(tx)$, respectively. Therefore, in this case, $s_m(x)$ can be
expressed as a sum of the terms $w(b,t)(c_b(tx)+d_b(tx))$, where $t|u^2$ and the coefficient
$w(b,t)=K_bt/u$.

Case 2:  Among the $k$ primes $u_1,u_2,\ldots,u_k$, there exists $q$ primes $u_{i_1},u_{i_2},\ldots,u_{i_q}$ with residue  $3$ modulo $4$.
To compute $s_m(x)$, we first consider the case when $q=1$, and then argue that
 the case $q>1$ can be dealt with in the same way.

Since $q=1$, we assume that the first $k-1$ odd primes
$u_{1},u_{2},\ldots,u_{k-1}$ satisfy $u_{i}\equiv 1\ (\bmod\  4)$
for $1\leq i\leq k-1$, and suppose that the last prime $u_{k}$ satisfies
$u_{k}\equiv 3\ (\bmod\  4)$, or
\[
\left(\frac{-b}{u_k}\right)=-\left(\frac{b}{u_k}\right).
\]

Similarly, we may define the indices $k_1$, $k_2$ and $k_3$  as
in Case 1 except that  $k_1+k_2+k_3=k-1$ since we now  have $k-1$ primes with residue  $1$ modulo $4$. This leads us to
consider two subcases according to whether $\left(\frac{-b}{u_k}\right)$ equals $0$. Keep in mind that $q=1$
in these two subcases.

On the one hand,   in order to use the two recursive relations (\ref{cbl}) and (\ref{dbl}) between $s_{m,n}$ and $s_{b,n}$,
we may assume that $\left(\frac{-b}{u_k}\right)=0$. Therefore, the term $u_k$ on the denominator and the same term
on the numerator cancel each other in (\ref{cbl}).
This argument also applies to the relation (\ref{dbl}). In other words,  there is no need to consider the occurrence of the
term $u_k^{2n+1}$ in (\ref{cbl}) and the terms $u_k^{2n}$ in (\ref{dbl}). In this sense, it remains to  consider the other $k-1$ primes
$u_{1},u_{2},\ldots,u_{k-1}$ such that $u_{i}\equiv 1\ (\bmod\  4)$
for $1\leq i\leq k-1$. By the argument in Case 1, we see again that  $s_m(x)$ can be
expressed as a sum of the terms $w(b,t)(c_b(tx)+d_b(tx))$, where $t|u^2$ and $w(b,t)=K_bt/u$.

On the other hand, we should consider the case for $\left(\frac{-b}{u_k}\right)=1$ or $\left(\frac{-b}{u_k}\right)=-1$. Since the proofs for these two situations are similar, we only give the proof of
 the case $\left(\frac{-b}{u_k}\right)=-1$. Assuming so,
  from  (\ref{cbl}) it follows that
\begin{align*}
c_{m,n}=&K_bu(u^{2})^{2n}\frac{\prod_{i=1}^{k_1}(u_{i}^{2n+1}-1)\prod_{j=1}^{k_2}(u_{k_1+j}^{2n+1}+1)}{\prod_{i=1}^{k_1+k_2}u_{i}^{2n+1}}\times \frac{u_k^{2n+1}+1}{u_k^{2n+1}}c_{b,n}\\
=&K_bu(u^{2})^{2n}
\left(\frac{\prod_{i=1}^{k_1}(u_{i}^{2n+1}-1)\prod_{j=1}^{k_2}(u_{k_1+j}^{2n+1}+1)}{\prod_{i=1}^{k_1+k_2}u_{i}^{2n+1}}+\frac{\prod_{i=1}^{k_1}(u_{i}^{2n+1}-1)\prod_{j=1}^{k_2}(u_{k_1+j}^{2n+1}+1)}{u_k^{2n+1}\prod_{i=1}^{k_1+k_2}u_{i}^{2n+1}}
\right) c_{b,n}.
\end{align*}
Let
\begin{equation}
g_b(x)=K_bu\sum\limits_{n\geq 0}
\left(\frac{\prod_{i=1}^{k_1}(u_{i}^{2n+1}-1)\prod_{j=1}^{k_2}(u_{k_1+j}^{2n+1}+1)}{\prod_{i=1}^{k_1+k_2}u_{i}^{2n+1}}\right)
c_{b,n}\frac{(u^{2}x)^{2n}}{(2n)!},\nonumber
\end{equation}
then we have
\begin{align}
c_m(x)=g_b(x)+\frac{1}{u_k}g_b(x/u_k)\label{gbc}.
\end{align}
In fact, by the argument in Case 1, we see that  $g_b(x)$ is a linear
combination of the terms $c_b(tx)$, where $t|u^2$  and the  coefficient of $c_b(tx)$ in linear combination is equal to  $K_bt/u$.

Similarly, by (\ref{dbl}) we find
\begin{align*}
d_{m,n}=&K_bu(u^{2})^{2n-1}\frac{\prod_{i=1}^{k_1}(u_{i}^{2n}-1)\prod_{j=1}^{k_2}(u_{k_1+j}^{2n}+1)}
{\prod_{i=1}^{k_1+k_2}u_{i}^{2n}}\times \frac{u_k^{2n}-1}{u_k^{2n}}d_{b,n}\\
=&K_bu(u^{2})^{2n-1}
\left(\frac{\prod_{i=1}^{k_1}(u_{i}^{2n}-1)\prod_{j=1}^{k_2}(u_{k_1+j}^{2n}+1)}
{\prod_{i=1}^{k_1+k_2}u_{i}^{2n}}-\frac{\prod_{i=1}^{k_1}(u_{i}^{2n}-1)\prod_{j=1}^{k_2}(u_{k_1+j}^{2n}+1)}
{u_k^{2n}\prod_{i=1}^{k_1+k_2}u_{i}^{2n}}
\right)d_{b,n}.
\end{align*}
Let
\begin{equation}
h_b(x)=K_bu\sum\limits_{n\geq 1}
\left(\frac{\prod_{i=1}^{k_1}(u_{i}^{2n}-1)\prod_{j=1}^{k_2}(u_{k_1+j}^{2n}+1)}{\prod_{i=1}^{k_1+k_2}u_{i}^{2n}}\right)
d_{b,n}\frac{(u^{2}x)^{2n-1}}{(2n-1)!}.\nonumber
\end{equation}
We get
\begin{align}
d_m(x)=h_b(x)-\frac{1}{u_k}h_b(x/u_k)\label{hbd}.
\end{align}
Again, from the reasoning in Case 1 it follows that $h_b(x)$ is a linear
combination of the terms $d_b(tx)$, where $t|u^2$  and the  coefficient of $d_b(tx)$ in linear combination is equal to  $K_bt/u$. Combining equations (\ref{gbc}) and (\ref{hbd}) yields the following formula
\begin{align*}
s_m(x)=g_b(x)+h_b(x)+\frac{1}{u_k}(g_b(x/u_k)-h_b(x/u_k)).
\end{align*}
Thus we see that
 $s_m(x)$ can also be expressed as a linear combination of the terms $w(b,t)(\pm c_b(tx)\pm d_b(tx)$, where $t|u^2$ and $w(b,t)=K_bt/u$.

Finally, as mentioned before we shall show that the justification for the above two subcases can be
applied to
 the case for  $q>1$. Let us give an example for  $q=2$.
Suppose that $u_{k-1}$ and $u_k$  are the last two primes such that $u_{k-1} \equiv 3\ (\bmod\  4)$,
$u_k \equiv 3\ (\bmod\  4)$, and $\left(\frac{-b}{u_{k-1}}\right)=-1$,
$\left(\frac{-b}{u_k}\right)=-1$. Then by using the same notation $g_b(x)$ and $h_b(x)$ as before, we can deduce the following formula for $s_m(x)$,
\begin{align*}
s_m(x)=&g_b(x)+h_b(x)+\frac{1}{u_{k-1}}(g_b(x/u_{k-1})-h_b(x/u_{k-1}))
+\frac{1}{u_k}(g_b(x/u_k)-h_b(x/u_k))\\
&+\frac{1}{u_{k-1}u_k}(g_b(x/u_{k-1}u_k)+h_b(x/u_{k-1}u_k)).
\end{align*}
Clearly, $s_m(x)$ can also be expressed as a linear combination of the terms $w(b,t)(\pm c_b(tx)\pm d_b(tx)$, where $t|u^2$ and $w(b,t)=K_bt/u$.
For other assumptions on $u_{k-1}$ and $u_k$, the computation can be similarly deduced.

In summary, $s_m(x)$ can  be expressed as a linear combination of the terms $w(b,t)(\pm c_b(tx)\pm d_b(tx)$,
where $t|u^2$ and $w(b,t)=K_bt/u$.
Since $b$ is square-free, $c_b(tx)$ and $d_b(tx)$ can be expressed as  a linear combination of the functions $\sec (btx)\cos((b-p)tx)$ and
 $\sec (btx)\sin(ptx)$, respectively, by using the formulas in the previous
section.   Hence the theorem holds.
\qed

Here we give three examples corresponding to the above three cases. For Case 1, suppose that $m=3(5\times 13)^2=3 (65)^2=3\times 4225=12675$.
Then we have
\begin{align*}
s_{12675}(x)=65s_3(4225x)-5s_3(325x)+13s_3(845x)-s_3(65x),
\end{align*}
where
\begin{equation}
s_3(x)=\sec (3x)(\sin 2x+\cos x).\nonumber
\end{equation}

For the first subcase of Case 2, suppose that  $m=6(5\times 3)^2=3 (15)^2=6\times 225=1350$. Then we get
\begin{align*}
s_{1350}(x)=15s_6(225x)-3s_6(45x),
\end{align*}
where
\begin{equation}
s_6(x)=\sec (
6x)( \cos 5x+\sin x)+\sec (6x)(\cos x+\sin 5x).\nonumber
\end{equation}

For the second subcase of Case 2, assume that $m=225=(5\times 3)^2$. We find
\begin{align*}
2S_{225}(x)=&15(c_1(225x)+d_1(225x))-3(c_1(45x)+d_1(45x))+5(c_1(75x)-d_1(75x)) \nonumber\\&-(c_1(15x)-d_1(15x)),
\end{align*}
where
\[
c_1(x)=\sec (x),\quad
d_1(x)=\tan (x).
\]

\section{Combinatorial Interpretation for  $s_{m,n}$}\label{sec4}

In this section, we aim to give a combinatorial interpretation for $s_{m,n}$ based on its generating function formula $s_m(x)$. Let us begin by recalling the known combinatorial interpretations of $s_{m,n}$ for $m=1,2,3,4$.

For $m=1$,  $(s_{1,n})_{n\geq 0}$ is called the sequence of Euler numbers. Let $E_n$ be the $n$-th Euler number, that is,
 the number of alternating permutations or up-down permutations of $[n]=\{1, 2, \ldots, n\}$, which are also
 called snakes of type $A_{n-1}$ by Arnol'd  \cite{Arno}. The following generating function is
  due to Andr\'e \cite{And}:
  \begin{equation}
  \sum\limits_{n\geq 0} E_n\frac{x^n}{n!} =  \sec x+\tan x,\nonumber
  \end{equation}
  Note that Springer also gave an explanation of the Euler numbers in terms
of the irreducible root system $A_{n-1}$ and derived the generating function of Andr\'e in this context.

For $m=2$, the sequence  $(s_{2,n})_{n\geq 0}$
   turns out to be the sequence of Springer numbers of the irreducible root system $B_n$ (\cite{Springer}). Purtill \cite{Purtill} has found an interpretation of this
 sequence. Let $P_n$ be the $n$-th entry of this sequence, whereas Purtill used the notation $E^{\pm}_n$. He has shown
        that $P_n$ equals the number of  Andr\'e signed permutations on $[n]$. On the other hand, it has been shown by Arnol'd  \cite{Arno} that $S_n$ also counts the
    number of  snakes of type $B_n$.  Hoffman \cite{Hoff} has derived the generating function of the number of snakes of type $B_n$ by giving a direct combinatorial proof.

For $m=3$, the sequence  $(s_{3,n})_{n\geq 0}$ has been studied by Ehrenborg and Readdy (\cite{Ehren}). Let $F_n$ denote the $n$-th Ehrenborg and Readdy number, which was denoted by $| ER_n|$, see
\cite{Hoff}.  It has been shown that $F_n$ equals the number $\Lambda$-alternating augmented $3$-signed permutations on $[n]$.
 Meanwhile, Hoffman \cite{Hoff} presented another combinatorial interpretation of the sequence in the case $m=3$ in terms of $ER_n$-snakes
 in the spirit of the snakes of type $A_{n-1}$ and $B_n$.

For $m=4$, the sequence  $(s_{4,n})_{n\geq 0}$ has also been studied by Ehrenborg and Readdy (\cite{Ehren2}).
They introduced the concept of non-augmented Andr\'e $4$-signed permutations on $[n]$ and proved that
such permutations are counted by $s_{4,n}$.

In principle, we have a combinatorial  interpretation of  $s_{m,n}$ for $m\geq 5$ based on the generating function $s_m(x)$.
Here are some definitions.
In \cite{Ehren}, Ehrenborg and Readdy defined a poset called the Sheffer poset, which can be viewed as a generalization of the binomial
poset introduced by Stanley \cite{Stan}. As an important example, they studied the $r$-cubical lattice, which is a set of ordered $r$-tuples
$(A_1,A_2,\ldots,A_r)$ of subsets from an infinite set $I$ together with the reverse inclusion order and a minimum element $\hat 0$ adjoined.
 Note that the $r$-cubical lattice has been studied by Metropolis, Rota, Strehl and White \cite{Rota}. Ehrenborg and Readdy further generalized
  the concept of $R$-labelings to linear edge-labelings. By considering the set of maximal chains in the interval $[\hat 0,\hat 1]$ on the Hasse diagram
  of the $r$-cubical
   lattice, they deduced a formula
for the number of $\Lambda$-alternating augmented $r$-signed permutations.

To introduce the definition of these permutations, Ehrenborg and Readdy constructed a linear edge-labeling on the Hasse diagram of
 the $r$-cubical lattice. To be more specific,
for an edge corresponding to the cover relation $A<B$ with $A\neq \hat 0$, let $(i,j)$ be its label where  $i$ equals the unique index such that
$A_i\neq B_i$ and $j$ takes the singleton element in $A_i-B_i$.   Let $G$ be the  label of the edge corresponding to  $\hat 0<A$,
whereas Ehrenborg and Readdy called it
 the special element. Then an augmented $r$-signed permutation is a list $(G,((i_1,j_1),(i_2,j_2),\ldots,(i_n,j_n)))$, where
  $i_1,i_2,\ldots,i_n\in [r]$ and $(j_1,j_2,\ldots,j_n)$ forms a permutation on $[n]$. In other words,  $r$-signed permutations are permutations
   on $[n]$ in which each
 element is assigned $r$ signs.

 To define the descent set of $r$-signed permutations, let $\Lambda$ be the set of such labels of the edges  on the Hasse diagram of
 the $r$-cubical lattice. It is easy to see that
 \[ \Lambda=([r]\times[n])\cup \{G\}.\]
   Let  $p$ be an integer such that $0\leq p\leq r$. For
 fixed $r$, $n$ and $p$, we can define a linear order   on $\Lambda$ which satisfies the following conditions
\begin{equation}
(i,j)<_{\Lambda}G\Rightarrow i\leq r-p,\label{psmaller}
\end{equation}
and
\begin{equation}
(i,j)>_{\Lambda}G\Rightarrow i> r-p,\label{pbigger}
\end{equation}
where  $(i,j)$ is the label of the edge corresponding to the cover relation $A<B$ such that $A$ covers
$\hat 0$, and $G$ is the special element.
For the remaining labels, we  may arrange them in the lexicographic order.  The descent set of an augmented $r$-signed permutation
$(G=g_0,g_1,\ldots,g_n)$ is defined as the set $\{k:g_{k-1}>_\Lambda g_k\}$,
where $g_k=(i_k,j_k)$ for $1\leq k\leq n$. Therefore, for an $\Lambda$-alternating augmented $r$-signed permutation,  that is, permutation having descent
set $\{2,4,\ldots\}$, it is necessary to have the condition  $G<(i_1,j_1)$,
or $i_1> r-p$. In other words, the labels above $G$ in this ordering are those whose first coordinate may take $p$ possible values from the set
$\{r-p+1,r-p+2,\ldots,r\}$.

If we denote the number of $\Lambda$-alternating
augmented $r$-signed permutations by $\Lambda_{r,p,n}$, then there are $\Lambda_{r,p,n}$ maximal chains with descent set $\{2,4,\ldots\}$ in the
interval $[\hat 0,A]$ with $A$ being an element of rank $n+1$. Based on this observation, Ehrenborg and Readdy have derived the following generating
function formula for $\Lambda_{r,p,n}$
 \begin{equation}
\Lambda_{r,p}(x)=\sum\limits_{n\geq 0}\Lambda_{r,p,n}\frac{x^n}{n!}=\frac{\cos ((r-p)x)+\sin (px)}{\cos (rx)}.\label{exponent}
\end{equation}

Now we interpret the generalized Euler and class numbers $s_{m,n}$ as follows.
When $m$ is square-free, it is easily seen that the formula for $s_{m,n}$ can be expressed as
 a linear combination of the number $\Lambda_{r,p_2,n}-\Lambda_{r,p_1,n}$ for $p_2>p_1$, and the number $\Lambda_{r,p_1,n}+\Lambda_{r,p_2,n}$. So we shall give combinatorial
 interpretations of these two numbers.

 For the first case, we consider the function
$\sec (rx)[\cos ((r-p_2)x)+\sin (p_2x)-\cos ((r-p_1)x)-\sin (p_1x)]$. It is the generating function for the number of $\Lambda$-alternating
augmented $r$-signed permutations $(G=g_0,(i_1,j_1),(i_2,j_2),\ldots,(i_n,j_n))$  with $r-p_2+1\leq i_1\leq r-p_1$, or the number of
maximal chains  in the interval $[\hat 0,\hat 1]$
whose first non-special edge has the label $(i_1,j_1)$ with $r-p_2+1\leq i_1\leq r-p_1$. Here  the first non-special edge is the one
 corresponding to the cover relation $A<B$ such that A covers $\hat 0$.
This gives a combinatorial interpretation of the number $\Lambda_{r,p_2,n}-\Lambda_{r,p_1,n}$.

 On the other hand, the numbers
 $\Lambda_{r,p_1,n}+\Lambda_{r,p_2,n}$ have the generating function
$\sec (rx)[\cos ((r-p_2)x)+\sin (p_2x)+\cos ((r-p_1)x)+\sin (p_1x)]$. Its combinatorial interpretation can be described in terms of  the $r$-cubical
lattice,  since the
 set of maximal
 chains in the interval $[\hat 0,\hat 1]$ constitutes a
  subposet of the $r$-cubical lattice. Let $P_1$ and $P_2$ denote these two subposets. Then $P_1$   consists of  maximal chains
  in $[\hat 0,\hat 1]$ whose first non-special edge is labeled with $(i_1,j_1)$, where  $r-p_1+1\leq i_1\leq r$. Similarly,  $P_2$ consists
   of the maximal chains in $[\hat 0,\hat 1]$  whose first
   non-special edge is labeled with   $(i_1,j_1)$, where
   $r-p_2+1\leq i_1\leq r$.
 Hence, $\Lambda_{r,p_1,n}+\Lambda_{r,p_2,n}$ equals  the number of
  maximal chains in the Hasse diagram of the
   disjoint union $P_1+P_2$.

 Therefore, the above two numbers can be endowed a combinatorial interpretation in terms of  $\Lambda$-alternating
augmented $r$-signed permutations, or
 in terms of the maximal chains in the Hasse diagram of the $r$-cubical
lattice.

 For the case when $m$ is not square-free, say, $m=bu^2$.
  We shall encounter the functions $\sec (bx)(\cos ((b-p)x)-\sin (px))$ in the formulas
for $s_{m,n}$. It is clear to see that these functions are the generating functions for the numbers $(-1)^n\Lambda_{b,p,n}$  after replacing
$x$ with $-x$ in  (\ref{exponent})
\begin{equation}
\Lambda_{b,p}(-x)=\sum\limits_{n\geq 0}(-1)^n\Lambda_{b,p,n}\frac{x^n}{n!}=\frac{\cos ((b-p)x)-\sin (px)}{\cos (bx)},\nonumber
\end{equation}
so that we may give an analogous interpretation of the numbers $(-1)^n\Lambda_{b,p,n}$.

\vspace{0.5cm}
 \noindent{\bf Acknowledgments.}  This work was supported by  the 973
Project, the PCSIRT Project of the Ministry of Education,  and the National Science
Foundation of China.

\end{document}